\input amstex
\input amsppt.sty
\magnification=\magstep1
\hsize=32truecc
\vsize=22.2truecm
\baselineskip=16truept
\NoBlackBoxes
\TagsOnRight \pageno=1 \nologo
\def\Z{\Bbb Z}
\def\N{\Bbb N}

\def\l{\left}
\def\r{\right}
\def\bg{\bigg}
\def\({\bg(}
\def\[{\bg\lfloor}
\def\){\bg)}
\def\]{\bg\rfloor}
\def\t{\text}
\def\f{\frac}

\def\ls{\leqslant}
\def\gs{\geqslant}

\def\nh{\noalign{\hrule}}
\def\hh{height4pt}

\def\om{&\omit &}

\def\Remark{\medskip\noindent{\it  Remark}}

\hbox {Ramanujan J. 42(2017), no.\,1, 59--67.}
\bigskip
\topmatter
\title A new theorem on the prime-counting function\endtitle
\author Zhi-Wei Sun\endauthor
\leftheadtext{Zhi-Wei Sun}
\affil Department of Mathematics, Nanjing University\\
 Nanjing 210093, People's Republic of China
  \\  zwsun\@nju.edu.cn
  \\ {\tt http://math.nju.edu.cn/$\sim$zwsun}
\endaffil
\abstract For $x>0$ let $\pi(x)$ denote the number of primes not exceeding $x$.
For integers $a$ and $m>0$, we determine when there is an integer $n>1$ with $\pi(n)=(n+a)/m$.
In particular, we show that for any integers $m>2$ and $a\ls\lceil e^{m-1}/(m-1)\rceil$ there is an integer $n>1$ with $\pi(n)=(n+a)/m$.
Consequently, for any integer $m>4$ there is a positive integer $n$ with $\pi(mn)=m+n$. We also pose several conjectures for further research; for example, we conjecture that
for each $m=1,2,3,\ldots$ there is a positive integer $n$ such that $m+n$ divides $p_m+p_n$, where $p_k$ denotes the $k$-th prime.
\endabstract
\thanks 2010 {\it Mathematics Subject Classification}. \,Primary 11A41, 11N05;
Secondary  05A15, 11A25, 11B39, 11B75.
\newline\indent Supported by the National Natural Science
Foundation of China.
\endthanks

\endtopmatter
\document

\heading{1. Introduction}\endheading

For $x>0$ let $\pi(x)$ denote the number of primes not exceeding $x$.
The function $\pi(x)$ is usually called the {\it prime-counting function}.
For $n\in\Z^+=\{1,2,3,\ldots\}$, let $p_n$ stand for the $n$-th prime.
By the Prime Number Theorem,
$$\pi(x)\sim\f x{\log x}\quad\t{as}\ x\to+\infty;$$
equivalently, $p_n\sim n\log n$ as $n\to+\infty$.
The asymptotic behaviors of $\pi(x)$ and $p_n$ have been intensively investigated by analytic number theorists.
Recently, the author [S15] formulated many conjectures on arithmetic properties of $\pi(x)$ and $p_n$ which depend on exact values of $\pi(x)$ or $p_n$.
For example, he conjectured that for any integer $n>1$, the number $\pi(kn)$ is prime for some $k=1,\ldots,n$.

In 1962, S. Golomb [G] found the following surprising property of $\pi(x)$: For any integer $k>1$ there is an integer $n>1$ with $n/\pi(n)=k$.
Along this line, we obtain the following general result.

\proclaim{Theorem 1.1} {\rm (i)} Let $m$ be any positive integer. For the set
$$S_m:=\l\{a\in\Z:\ \pi(n)=\f{n+a}m\ \t{for some integer}\ n>1\r\},\tag1.1$$
we have
$$S_m=\{\ldots,-2,-1,\ldots,S(m)\},\tag1.2$$
where
$$S(m):=\max\{km-p_k:\ k\in\Z^+\}=\max\{km-p_k:\ k=1,2,\ldots,\lfloor e^{m+1}\rfloor\}.\tag1.3$$

{\rm (ii)} We have
$$(m-1)S(m+1)>mS(m)\quad\t{for any}\ m\in\Z^+.\tag1.4$$
Also,
$$\f{e^{m-1}}{m-1}<S(m)<(m-1)e^{m+1}\quad\t{for all}\ m=3,4,\ldots,\tag1.5$$
and hence
$$\lim_{m\to+\infty}\root{m}\of{S(m)}=e.\tag1.6$$
\endproclaim
\Remark\ 1.1. For any integer $m\gs2$, we have $S(m)\gs m-p_1\gs0$ and hence Theorem 1.1 yields Golomb's result $0\in S_m$.
In view of (1.5), for each $m=3,4,\ldots$, the least $k\in\Z^+$ with $km-p_k=S(m)$ is greater than $e^{m-1}/(m-1)^2$.

\medskip

\proclaim{Corollary 1.1} Let $m>0$ and $a\ls m^2-m-1$ be integers. Then there is an integer $n>1$
with $\pi(n)=(n+a)/m$, i.e.,
$$\pi(mn-a)=n\quad\t{for some}\ n\in\Z^+.\tag1.7$$
\endproclaim
\Remark\ 1.2. For any positive integer $m$,  if we let $n$ be the number of primes not exceeding the $m$-th composite number, then $\pi(m+n)=n$.

\proclaim{Corollary 1.2} For any integer $m>4$, there is a positive integer $n$ such that
$$\pi(mn)=m+n.\tag1.8$$
\endproclaim
\Remark\ 1.3. Let $n$ be any positive integer. Clearly $\pi(n)<n+1$ and $\pi(2n)\ls n<n+2$. Observe that
$$\align 2n=&\l\lfloor\f{3n}2\r\rfloor+n-\l\lfloor \f n2\r\rfloor
\\=&|\{1\ls k\ls 3n:\ 2\mid k\}|+|\{1\ls k\ls 3n:\ 3\mid k\}|-|\{1\ls k\ls 3n:\ 6\mid k\}|
\\=&|\{1\ls k\ls 3n:\ \gcd(k,6)>1\}|
\\\ls&|\{1\ls k\ls 3n:\ k \ \t{is not prime}\}|+1=3n-\pi(3n)+1
\endalign$$
and hence $\pi(3n)\ls n+1<n+3$. As $k:=\pi(4n)\gs2$, we have $4n\gs p_k\gs k(\log k+\log\log k-1)$ by [D]. If $n\gs45$, then
$\log k+\log\log k\gs 5$ and hence $\pi(4n)=k\ls n<n+4$. We can easily verify that $\pi(4n)<n+4$ if $n\ls 44$.
\medskip

Recall that the well-known Fibonacci numbers $F_n\ (n\in\N=\{0,1,2,\ldots\})$ are given by
$$F_0=0,\ F_1=1,\ \t{and}\ F_{k+1}=F_k+F_{k-1}\ (k=1,2,3,\ldots).$$

\proclaim{Corollary 1.3} For any integer $m>3$, there is a positive integer $n$ such that
$$\pi(mn)=F_m+n.\tag1.9$$
\endproclaim

A positive integer $n$ is called a {\it practical number} if every $m=1,\ldots,n$ can be expressed as a sum of some distinct (positive) divisors of $n$.
The only odd practical number is $1$. The distribution of practical numbers is quite similar to that of prime numbers.
For $x>0$ let $P(x)$ denote the number of practical numbers not exceeding $x$. Similar to the Prime Number Theorem, we have
$$P(x)\sim c\f x{\log x}\quad\t{for some constant}\ c>0,$$
which was conjectured by M. Margenstern [M] in 1991 and proved by A. Weingartner [W] in 2014.
In view of this, our method to prove Theorem 1.1(i) allows us to deduce for any positive integer $m$ the equality
$$\l\{a\in\Z:\ P(n)=\f{n+a}m\quad\t{for some}\ n\in\Z^+\r\}=\{\ldots,-2,-1,0,\ldots,T(m)\},\tag1.10$$
where $T(m)=\max\{km-q_k:\ k\in\Z^+\}$ with $q_k$ the $k$-th practical number.

We are going to show Theorem 1.1 in the next section. Section 3 contains our proofs of Corollaries 1.1-1.3 and related numerical tables.
In Section 4 we pose several conjectures for further research.

\heading{2. Proof of Theorem 1.1}\endheading

\medskip
\noindent{\it Proof of Theorem} 1.1(i). By [D],
$$p_k\gs k(\log k+\log\log k-1)\quad\t{for all}\ k=2,3,\ldots.$$
So, for any integer $k>e^{m+1}$, we have
$$km-p_k\ls k(m-\log(k\log k)+1)<0\quad\t{and hence}\ km-p_k\ls -1\ls m-p_1.$$
Therefore $S(m)=\max\{km-p_k:\ k=1,2,\ldots,\lfloor e^{m+1}\rfloor\}.$

For any $a\in S_m$, there is an integer $n>1$ such that $k:=\pi(n)=(n+a)/m$ and hence $a=km-n\ls km-p_k\ls S(m)$.

Define $I_k:=\{km-p_{k+1}+1,\ldots,km-p_{k+1}+m\}$ for all $k\in\N$.
As $\min I_0=-1$, and $\min I_{k+1}\ls\max I_k$ for all $k\in\N$, we see that
$\bigcup_{k\in\N}I_k\supseteq\{-1,\ldots,S(m)\}$.
Note that $\max I_k\ls S(m)$ and $km-p_{k+1}\to-\infty$. If $a$ is an integer with $\max I_{k+1}<a<\min I_k$, then for $n=(k+1)m-a$ we have
$p_{k+1}+m-1<n<p_{k+2}-m$, hence $a\in S_m$ since $\pi(n)=k+1=(n+a)/m$. Therefore
$$\(\bigcup_{k\in\N}I_k\)\cup S_m=\{\ldots,-2,-1,\ldots,S(m)\}.\tag2.1$$

Now suppose that $a$ is an integer with $a\ls S(m)$ and $a\not\in S_m$. We want to deduce a contradiction. In light of (2.1), for some $k\in\N$ we have
$$a\in I_k=\{km-p_{k+1}+1,\ldots,km-p_{k+1}+m\}.\tag2.2$$
Write $a=km+r$ with $1-p_{k+1}\ls r\ls m-p_{k+1}$. We claim that
$$\f{n+r}{\pi(n)-k}=m\quad\t{for some integer}\ n\gs p_{k+1}.\tag2.3$$
This is obvious for $m=p_{k+1}+r$ since
$$\f{p_{k+1}+r}{\pi(p_{k+1})-k}=p_{k+1}+r.$$
Below we assume $m>p_{k+1}+r$. As $\pi(n)\sim n/\log n$, we see that
$$\lim_{n\to+\infty}\f{n+r}{\pi(n)-k}=+\infty.$$
So, we may choose the least integer $n\gs p_{k+1}$ with $(n+r)/(\pi(n)-k)\gs m$.
Clearly $n\not=p_{k+1}$, thus $n-1\gs p_{k+1}$ and hence
$$\f{n+r}{\pi(n)-k}\gs m>\f{(n-1)+r}{\pi(n-1)-k}\tag2.4$$
by the choice of $n$. Set
$$s=n-1+r\quad\t{and}\quad t=\pi(n-1)-k.$$
As $n-1\gs p_{k+1}$, we have $t\gs 1$. Note also that
$$\align s-t=&n-1+r-(\pi(n-1)-k)
\\\gs& n-1+(1-p_{k+1})-\pi(n-1)+k
\\=&(n-1-p_{k+1})-(\pi(n-1)-\pi(p_{k+1}))
\\=&|\{p_{k+1}<d\ls n-1:\ d\ \t{is composite}\}|
\\\gs&0.
\endalign$$
If $n$ is prime, then $\pi(n)=\pi(n-1)+1$ and hence
$$\f{n+r}{\pi(n)-k}=\f{s+1}{t+1}\ls \f st=\f{n-1+r}{\pi(n-1)-k}$$
which contradicts (2.4). Thus $n$ is not prime and hence
$$n+r\gs m(\pi(n)-k)= m(\pi(n-1)-k)>n-1+r.$$
It follows that
$$\f{n+r}{\pi(n)-k}=m.$$

By the claim (2.3), for some integer $n\gs p_{k+1}$ we have
$$\pi(n)=k+\f{n+r}m=\f{n+a}m.$$
Therefore $a\in S_m$, which contradicts the supposition.

In view of the above, we have completed the proof of Theorem 1.1(i). \qed

\medskip
\noindent{\it Proof of Theorem} 1.1(ii). For any given $m\in\Z^+$, we may choose $k\in\Z^+$ with $km-p_k=S(m)$, and hence
$$\align(m-1)S(m+1)\gs &(m-1)(k(m+1)-p_k)=(m-1)S(m)+k(m-1)
\\>&(m-1)S(m)+km-p_k=mS(m).\endalign$$
This proves (1.4).

Clearly (1.6) follows from (1.5). Let $m>2$ be an integer. As $p_k>k$ for $k\in\Z^+$, we have
$S(m)<(m-1)e^{m+1}$ by (1.3). So it remains to show $j:=\lfloor e^{m-1}/(m-1)\rfloor<S(m)$.

For $m=3$, we clearly have $j=3<3\times3-p_3\ls S(3)$.

Below we assume $m\gs 4$. Then $j\gs6$ and hence
$$p_j\ls j(\log j+\log\log j)$$
by [RS, (3.13)] and [D, Lemma 1]. Clearly
$$\log j\ls \log\f{e^{m-1}}{m-1}=m-1-\log(m-1)<m-1,$$
and thus
$$jm-p_j\gs j(m-\log j)-j\log\log j>j(1+\log(m-1))-j\log(m-1)=j.$$
Therefore $j<S(m)$ as desired. \qed

\heading{3. Proofs of Corollaries 1.1-1.3 and related data}\endheading

\medskip
\noindent{\it Proof of Corollary} 1.1. By Theorem 1.1, it suffices to show that $m^2-m-1\ls S(m)$.

For $m\ls 5$, we have $m^2-m-1\ls k m-p_k$ for some $k\in\Z^+$. In fact,
$$\gather1^2-1-1=1\times1-p_1,
\ 2^2-2-1=2\times2-p_2,\ 3^2-3-1=5=4\times3-p_4,
\\ 4^2-4-1=11=6\times4-p_6\ \ \t{and}\ \ 5^2-5-1=19<8\times5-p_8=21.
\endgather$$
For $m\gs6$, we have $m^2-m-1<e^{m-1}/(m-1)$ and hence $m^2-m-1<S(m)$ by (1.5).
This concludes the proof. \qed
\medskip

As $S(m)=\max\{km-p_k:\ k=1,\ldots,\lfloor e^{m+1}\rfloor\}$, we can determine the exact values of $S(m)$ for small positive integers $m$.

\medskip
\centerline{Values of $S(m)$ for $m=1,\ldots,17$}
\smallskip
\centerline{\vbox{\offinterlineskip
\halign{\vrule#&\ \ #\ \hfill   &&\vrule#&\ \ \hfill#\ \
\cr\nh\cr \hh\om\om\om\om\om\om\om\om\om\om\om
\cr &$m$& &$1$& &$2$& &$3$& &$4$& &$5$&  &$6$& &$7$& &$8$& &$9$& &$10$&
\cr \hh\om\om\om\om\om\om\om\om\om\om\om
\cr\nh\cr \hh\om\om\om\om\om\om\om\om\om\om\om
\cr &$S(m)$& &$-1$& &$1$& &$5$& &$13$& &$37$& &$83$& &$194$& &$469$& &$1111$& &$2743$&
\cr \hh\om\om\om\om\om\om\om\om\om\om\om
\cr \nh\cr }}}
\medskip

\centerline{\vbox{\offinterlineskip
\halign{\vrule#&\ \ #\ \hfill   &&\vrule#&\ \ \hfill#\ \
\cr\nh\cr \hh\om\om\om\om\om\om\om\om
\cr &$m$& &$11$& &$12$& &$13$& &$14$&  &$15$& &$16$& &$17$&
\cr \hh\om\om\om\om\om\om\om\om
\cr\nh\cr \hh\om\om\om\om\om\om\om\om
\cr &$S(m)$& &$6698$& &$16379$& &$40543$& &$101251$& &$254053$& &$640483$& &$1622840$&
\cr \hh\om\om\om\om\om\om\om\om
\cr \nh\cr }}}
\medskip

In the following table, for each $m=2,\ldots,20$ we give the least integer $n>1$ with $\pi(n)=(n-1)/m$ as well as the least integer $n>1$ with $\pi(n)=(n+m-1)/m$.

\bigskip
\centerline{}
\smallskip
\centerline{\vbox{\offinterlineskip
\halign{\vrule#&\ \ #\ \hfill   &&\vrule#&\ \ \hfill#\ \
\cr\nh\cr \hh\om\om\om
\cr &$m$& &$\t{Least}\ n>1\ \t{with}\ \pi(n)=\f{n-1}m$& &$\t{Least}\ n>1\ \t{with}\ \pi(n)=\f{n+m-1}m$&
\cr \hh\om\om\om
\cr\nh\cr \hh\om\om\om
\cr &$2$& &$9$& &$3$&
\cr \hh\om\om\om
\cr\nh\cr \hh\om\om\om
\cr &$3$& &$28$& &$4$&
\cr \hh\om\om\om
\cr\nh\cr \hh\om\om\om
\cr &$4$& &$121$& &$93$&
\cr \hh\om\om\om
\cr\nh\cr \hh\om\om\om
\cr &$5$& &$336$& &$306$&
\cr \hh\om\om\om
\cr\nh\cr \hh\om\om\om
\cr &$6$& &$1081$& &$1003$&
\cr \hh\om\om\om
\cr\nh\cr \hh\om\om\om
\cr &$7$& &$3060$& &$2997$&
\cr \hh\om\om\om
\cr\nh\cr \hh\om\om\om
\cr &$8$& &$8409$& &$8361$&
\cr \hh\om\om\om
\cr\nh\cr \hh\om\om\om
\cr &$9$& &$23527$& &$23518$&
\cr \hh\om\om\om
\cr\nh\cr \hh\om\om\om
\cr &$10$& &$64541$& &$64531$&
\cr \hh\om\om\om
\cr\nh\cr \hh\om\om\om
\cr &$11$& &$175198$& &$175187$&
\cr \hh\om\om\om
\cr\nh\cr \hh\om\om\om
\cr &$12$& &$480865$& &$480817$&
\cr \hh\om\om\om
\cr\nh\cr \hh\om\om\om
\cr &$13$& &$1304499$& &$1303004$&
\cr \hh\om\om\om
\cr\nh\cr \hh\om\om\om
\cr &$14$& &$3523885$& &$3523871$&
\cr \hh\om\om\om
\cr\nh\cr \hh\om\om\om
\cr &$15$& &$9557956$& &$9557746$&
\cr \hh\om\om\om
\cr\nh\cr \hh\om\om\om
\cr &$16$& &$25874753$& &$25874737$&
\cr \hh\om\om\om
\cr\nh\cr \hh\om\om\om
\cr &$17$& &$70115413$& &$70115311$&
\cr \hh\om\om\om
\cr\nh\cr \hh\om\om\om
\cr &$18$& &$189961183$& &$189961075$&
\cr \hh\om\om\om
\cr\nh\cr \hh\om\om\om
\cr &$19$& &$514272412$& &$514272393$&
\cr \hh\om\om\om
\cr\nh\cr \hh\om\om\om
\cr &$20$& &$1394193581$& &$1394193361$&
\cr \hh\om\om\om
\cr \nh\cr }}}
\smallskip
\bigskip

\noindent{\it Proof of Corollary} 1.2. Note that $\pi(5\times9)=5+9$ and $\pi(6\times7)=6+7$.

Now we assume $m\gs7$. Then $m^2<e^{m-1}/(m-1)$. By Theorem 1.1, there is a positive integer $N$ with $\pi(N)=(N+m^2)/m$.
Clearly $n=N/m\in\Z^+$ and $\pi(mn)=(mn+m^2)/m=m+n$. This concludes the proof. \qed

\medskip
\centerline{Smallest $n=s(m)$ with $\pi(mn)=m+n$ for $5\ls m\ls 21$}
\smallskip
\centerline{\vbox{\offinterlineskip
\halign{\vrule#&\ \ #\ \hfill   &&\vrule#&\ \ \hfill#\ \
\cr\nh\cr \hh\om\om\om\om\om\om\om\om
\cr &$m$& &$5$& &$6$& &$7$& &$8$& &$9\sim14$&  &$15$& &$16$&
\cr \hh\om\om\om\om\om\om\om\om
\cr\nh\cr \hh\om\om\om\om\om\om\om\om
\cr &$s(m)$& &$9$& &$7$& &$6$& &$998$& &$5$& &$636787$& &$1617099$&
\cr \hh\om\om\om\om\om\om\om\om
\cr \nh\cr }}}
\medskip

\centerline{\vbox{\offinterlineskip
\halign{\vrule#&\ \ #\ \hfill   &&\vrule#&\ \ \hfill#\ \
\cr\nh\cr \hh\om\om\om\om\om\om
\cr &$m$& &$17$& &$18$& &$19$& &$20$& &$21$&
\cr \hh\om\om\om\om\om\om
\cr\nh\cr \hh\om\om\om\om\om\om
\cr &$s(m)$& &$4124188$&  &$10553076$& &$5$& &$5$& &$179992154$&
\cr \hh\om\om\om\om\om\om
\cr \nh\cr }}}
\smallskip
\bigskip

\noindent{\it Proof of Corollary} 1.3. Observe that
$$\gather\pi(4\times5)=F_4+5,\ \pi(5\times9)=F_5+9,\ \pi(6\times12)=F_6+12,
\\ \pi(7\times16)=F_7+16\ \ \t{and}\ \ \pi(8\times25)=F_8+25.
\endgather$$

Now we assume $m\gs9$. Then $mF_m<e^{m-1}/(m-1)$. By Theorem 1.1, there is a positive integer $N$ with $\pi(N)=(N+mF_m)/m$.
Note that $n=N/m\in\Z^+$ and $\pi(mn)=(mn+mF_m)/m=F_m+n$. This concludes the proof. \qed

\medskip
\centerline{Least $n=f(m)$ with $\pi(mn)=F_m+n$ for $4\ls m\ls 22$}
\smallskip
\centerline{\vbox{\offinterlineskip
\halign{\vrule#&\ \ #\ \hfill   &&\vrule#&\ \ \hfill#\ \
\cr\nh\cr \hh\om\om\om\om\om\om\om\om\om\om\om\om
\cr &$m$& &$4$& &$5$& &$6$& &$7$& &$8$&  &$9$& &$10$& &$11$& &$12$& &$13$& &$14$&
\cr \hh\om\om\om\om\om\om\om\om\om\om\om\om
\cr\nh\cr \hh\om\om\om\om\om\om\om\om\om\om\om\om
\cr &$f(m)$& &$5$& &$9$& &$12$& &$16$& &$25$& &$45$& &$68$& &$116$& &$183$& &$287$& &$457$&
\cr \hh\om\om\om\om\om\om\om\om\om\om\om\om
\cr \nh\cr }}}
\medskip

\centerline{\vbox{\offinterlineskip
\halign{\vrule#&\ \ #\ \hfill   &&\vrule#&\ \ \hfill#\ \
\cr\nh\cr \hh\om\om\om\om\om\om\om\om\om
\cr &$m$& &$15$& &$16$& &$17$& &$18$& &$19$&  &$20$& &$21$& &$22$&
\cr \hh\om\om\om\om\om\om\om\om\om
\cr\nh\cr \hh\om\om\om\om\om\om\om\om\om
\cr &$f(m)$& &$628346$& &$1600659$& &$1942$& &$3133$& &$5028$& &$8131$& &$13100$& &$21142$&
\cr \hh\om\om\om\om\om\om\om\om\om
\cr \nh\cr }}}
\bigskip

\heading{4. Some conjectures}\endheading

In view of Theorem 1.1, we pose the following conjecture.

\proclaim{Conjecture 4.1} {\rm (i)} Let $m$ be any positive integer. Then $p_k-km$ is prime for some $k\in\Z^+$,
and $p_k-km$ is a square for some $k\in\Z^+$. If $m>2$, then $km-p_k$ is prime for some $k\in\Z^+$, and $km-p_k$ is a square for some $k\in\Z^+$.

{\rm (ii)} The sequence $\root m\of{S(m)}\ (m=1,2,3,\ldots)$ is strictly increasing.
\endproclaim
\Remark\ 4.1. See [S14, A247278, A247893 and A247895] for some sequences related to part (i); for example, $29\times 5-p_{29}=145-109=6^2$ and $p_{12}-12\times 3=1^2$.
The second part of Conjecture 4.1 arises naturally in the spirit of [S13].
\medskip

Golomb's result [G] indicates that for any integer $m\gs2$ we have $\pi(mn)=n\,(=mn/m)$ for some $n\in\Z^+$. Motivated by this and Corollary 1.2, we pose the following conjecture
related to Euler's totient function $\varphi$.

\proclaim{Conjecture 4.2} Let $m$ be any positive integer. Then $\pi(mn)=\varphi(n)$ for some $n\in\Z^+$.
Also, $\pi(mn)=\varphi(m)+\varphi(n)$ for some $n\in\Z^+$, and $\pi(mn)=\varphi(m+n)$ for some $n\in\Z^+$.
\endproclaim
\Remark\ 4.2. Our method to establish Theorem 1.1 does not work for this conjecture.

\medskip
\centerline{Least $n\in\Z^+$ with $\pi(mn)=\varphi(n)$ for $m\ls 18$}
\smallskip
\centerline{\vbox{\offinterlineskip
\halign{\vrule#&\ \ #\ \hfill   &&\vrule#&\ \ \hfill#\ \
\cr\nh\cr \hh\om\om\om\om\om\om\om\om\om\om\om\om
\cr &$m$& &$1$& &$2$& &$3$& &$4$& &$5$&  &$6$& &$7$& &$8$& &$9$& &$10$& &$11$&
\cr \hh\om\om\om\om\om\om\om\om\om\om\om\om
\cr\nh\cr \hh\om\om\om\om\om\om\om\om\om\om\om\om
\cr &$n$& &$2$& &$1$& &$13$& &$31$& &$73$& &$181$& &$443$& &$2249$& &$238839$& &$6473$& &$3001$&
\cr \hh\om\om\om\om\om\om\om\om\om\om\om\om
\cr \nh\cr }}}
\medskip

\centerline{\vbox{\offinterlineskip
\halign{\vrule#&\ \ #\ \hfill   &&\vrule#&\ \ \hfill#\ \
\cr\nh\cr \hh\om\om\om\om\om\om\om\om
\cr &$m$& &$12$& &$13$& &$14$& &$15$& &$16$&  &$17$& &$18$&
\cr \hh\om\om\om\om\om\om\om\om
\cr\nh\cr \hh\om\om\om\om\om\om\om\om
\cr &$n$& &$40123$& &$108539$& &$251707$& &$637321$& &$7554079$& &$4124437$& &$241895689$&
\cr \hh\om\om\om\om\om\om\om\om
\cr \nh\cr }}}
\bigskip

\newpage
\centerline{Least $n\in\Z^+$ with $\pi(mn)=\varphi(m)+\varphi(n)$ for $m\ls 18$}
\smallskip
\centerline{\vbox{\offinterlineskip
\halign{\vrule#&\ \ #\ \hfill   &&\vrule#&\ \ \hfill#\ \
\cr\nh\cr \hh\om\om\om\om\om\om\om\om\om\om\om\om
\cr &$m$& &$1$& &$2$& &$3$& &$4$& &$5$&  &$6$& &$7$& &$8$& &$9$& &$10$& &$11$&
\cr \hh\om\om\om\om\om\om\om\om\om\om\om\om
\cr\nh\cr \hh\om\om\om\om\om\om\om\om\om\om\om\om
\cr &$n$& &$6$& &$2$& &$2$& &$23$& &$3$& &$1$& &$3$& &$1033$& &$2$& &$6449$& &$15887$&
\cr \hh\om\om\om\om\om\om\om\om\om\om\om\om
\cr \nh\cr }}}
\medskip

\centerline{\vbox{\offinterlineskip
\halign{\vrule#&\ \ #\ \hfill   &&\vrule#&\ \ \hfill#\ \
\cr\nh\cr \hh\om\om\om\om\om\om\om\om
\cr &$m$& &$12$& &$13$& &$14$& &$15$& &$16$&  &$17$& &$18$&
\cr \hh\om\om\om\om\om\om\om\om
\cr\nh\cr \hh\om\om\om\om\om\om\om\om
\cr &$n$& &$1$& &$100169$& &$268393$& &$636917$& &$2113589$& &$70324093$& &$1$&
\cr \hh\om\om\om\om\om\om\om\om
\cr \nh\cr }}}
\bigskip

\medskip
\centerline{Least $n\in\Z^+$ with $\pi(mn)=\varphi(m+n)$ for $m\ls 20$}
\smallskip
\centerline{\vbox{\offinterlineskip
\halign{\vrule#&\ \ #\ \hfill   &&\vrule#&\ \ \hfill#\ \
\cr\nh\cr \hh\om\om\om\om\om\om\om\om\om\om\om\om\om\om
\cr &$m$& &$1$& &$2$& &$3$& &$4$& &$5$&  &$6$& &$7$& &$8$& &$9$& &$10$& &$11$& &$12$& &$13$&
\cr \hh\om\om\om\om\om\om\om\om\om\om\om\om\om\om
\cr\nh\cr \hh\om\om\om\om\om\om\om\om\om\om\om\om\om\om
\cr &$n$& &$3$& &$2$& &$1$& &$91$& &$6$& &$5$& &$1$& &$5$& &$1$& &$8041$& &$15870$& &$39865$& &$1$&
\cr \hh\om\om\om\om\om\om\om\om\om\om\om\om\om\om
\cr \nh\cr }}}
\medskip

\centerline{\vbox{\offinterlineskip
\halign{\vrule#&\ \ #\ \hfill   &&\vrule#&\ \ \hfill#\ \
\cr\nh\cr \hh\om\om\om\om\om\om\om\om
\cr &$m$&  &$14$& &$15$& &$16$&  &$17$& &$18$& &$19$& &$20$&
\cr \hh\om\om\om\om\om\om\om\om
\cr\nh\cr \hh\om\om\om\om\om\om\om\om
\cr &$n$&  &$251625$& &$637064$& &$1829661$& &$4124240$& &$10553093$& &$1$& &$69709253$&
\cr \hh\om\om\om\om\om\om\om\om
\cr \nh\cr }}}
\bigskip

For $n\in\Z^+$ let $\sigma(n)$ denote the number of (positive) divisors of $n$. We also formulate the following conjecture motivated by Conjecture 4.2.

\proclaim{Conjecture 4.3} For any integer $m>1$, there is a positive integer $n$ with $\pi(mn)=\sigma(n)$. Also, for any integer $m>4$, $\pi(mn)=\sigma(m)+\sigma(n)$
for some $n\in\Z^+$, and $\pi(mn)=\sigma(m+n)$ for some $n\in\Z^+$.
\endproclaim

{\it Example} 4.1. The least $n\in\Z^+$ with $\pi(23n)=\sigma(n)$ is $8131355$, the least $n\in\Z^+$ with $\pi(39n)=\sigma(39)+\sigma(n)$
is $75999272$, and the least $n\in\Z^+$ with $\pi(30n)=\sigma(30+n)$ is $39298437$.
\medskip

Now we pose one more conjecture which is motivated by Corollary 1.2.

\proclaim{Conjecture 4.4} Let $m$ be any positive integer. Then $m+n$ divides $p_m+p_n$ for some $n\in\Z^+$. Moreover, we may require $n<m(m-1)$ if $m>2$.
\endproclaim
\Remark\ 4.3. We have verified this for all $m=1,\ldots,10^5$, see [S14, A247824] for related data. We also conjecture that for any $m\in\Z^+$
there is a positive integer $n$ such that $\pi(mn)$ divides $p_m+p_n$, see [S14, A247793] for related data.
\medskip

{\it Example} 4.2. The least $n\in\Z^+$ with $2+n$ dividing $p_2+p_n$ is $5$. For $m=79276$, the least $n\in\Z^+$ with $m+n$ dividing $p_m+p_n$
is $3141281384>3\times10^9$.

\medskip



\widestnumber\key{S14}

 \Refs


\ref\key D\by P. Dusart\paper The $k$th prime is greater than $k(\log k+\log\log k-1)$ for $k\gs2$
\jour Math. Comp. \vol 68\yr 1999\pages 411--415\endref

\ref\key G\by S. W. Golomb\paper On the ratio of $N$ to $\pi(N)$\jour Amer. Math. Monthly\vol 69\yr 1962\pages 36--37\endref

\ref\key M\by M. Margenstern\paper Les nombres pratiques: th\'eorie, observations et conjectures\jour J. Number Theory\vol 37\yr 1991\pages 1--36\endref

\ref\key RS\by J. B. Rosser and L. Schoenfeld\paper Approximate formulas for some functions of prime numbers
\jour Illinois J. Math. \vol 6\yr 1962\pages 64--94\endref

\ref\key S13\by Z.-W. Sun\paper Conjectures involving arithmetical sequences\jour in: S. Kanemitsu, H. Li and J. Liu (eds.), Number Theory: Arithmetic in Shangri-La,
Proc. 6th China-Japan Seminar (Shanghai, August 15-17, 2011), World Sci., Singapore, 2013, pp. 244-258\endref

\ref\key S14\by Z.-W. Sun\paper {\rm Sequences  A247278, A247793, A247824, A247893 and A247895 in OEIS (On-Line Encyclopedia of Integer Sequences)}
\jour {\tt http://oeis.org}\endref

\ref\key S15\by Z.-W. Sun\paper Problems on combinatorial properties of primes\jour in: M. Kaneko, S. Kanemitsu and J. Liu (eds.), Plowing and Starring through High Wave Forms, Proc. 7th China-Japan Seminar on Number Theory
(Fukuoka, Oct. 28--Nov. 1, 2013), Ser. Number Theory Appl., Vol. 11, World Sci., Singapore, 2015, pp. 169--187\endref

\ref\key W\by A. Weingartner\paper Practical numbers and the distribution of divisors\jour Q. J. Math. \vol 66\yr 2015\pages 743--758\endref

\endRefs
\enddocument